\documentclass[a4paper,12pt]{amsart}
\usepackage{amsfonts}
\def\T{{\mathcal{T}}}
\def\M{{\mathcal{M}}}
\def\O{{\mathcal{O}}}
\def\V{{\mathcal{V}}}

\def\oM{{{\overline{\mathcal{M}}}}}

\def\CP1{{\mathbb{C}\mathrm{P}^1}}

\def\l{\langle}
\def\r{\rangle}

\newtheorem{theorem}{Theorem}
\newtheorem{lemma}{Lemma}

\title{Intersections in genus 3 and the Boussinesq hierarchy}

\author{S.~V.~Shadrin}

\thanks{Partially supported by the grants RFBR 01-01-00660 and RFBR 02-01-22004.}

\begin{document}

\begin{abstract}
In this paper we prove that the enlarged Witten's conjecture is true in
the case of the Boussinesq hierarchy for correlators
$\l\tau_{n,m}\tau_{0,1}^k\tau_{0,0}^l\r_3$.
\end{abstract}

\maketitle





\section{Introduction}

\subsection{}
In this paper we study a very special case of the enlarged Witten's conjecture.
In~\cite{w0,w} Witten conjectured that the generating functions of certain intersection
numbers $\l\tau_{n_1,m_1}\dots\tau_{n_s,m_s}\r_g$ on the moduli spaces of curves are the
string solutions to the $n$-KdV hierarchies.

This conjecture was completely proved in the case of $2$-KdV hierarchy in~\cite{k}
and~\cite{op} (see also~\cite{l} for some comments). For the $n$-KdV, $n=3,4,\dots$,
this conjecture was proved only in genus zero ($g=0$), see~\cite{w}.

As for the higher genera, the situation looks like follows. In genera $1$ and $2$, there
exist topological recursion relations (see~\cite{g}) allowing to calculate intersection
numbers participating in the Witten's conjecture (that's why we consider genus $3$ in this
paper). For an arbitrary genus there exist an algorithm (see~\cite{s}) allowing to calculate
(in a sense) only the intersection numbers like
$\l\tau_{n,m}\tau_{0,m_1}\dots\tau_{0,m_s}\r_g$ (with descendants only at one point).

\subsection{}
In this paper we restrict ourself to the case of $3$-KdV hierarchy. We consider
intersections in genus $3$ with descendants only at one point. For these numbers we
give some relations slightly generalizing relations from~\cite{s}. Using these relation
we prove an expression of intersection numbers in genus $3$ via intersection numbers in
genus $0$. The same expression is predicted by the $3$-KdV (Boussinesq) hierarchy.
Since the conjecture is already prove in genuz zero, our theorem (theorem~\ref{theorem3})
follows that the conjecture is true for the intersection numbers
$\l\tau_{n,m}\tau_{0,1}^k\tau_{0,0}^l\r_3$ (of course, in the case of the
$3$-KdV hierarchy).

The paper is organized as follows. In section~\ref{witcon} we recall the Witten's
conjecture (the particular case we need here). In section~\ref{relint} we define some
auxiliary intersection numbers and give some relations. In section~\ref{calcul} we
outline the calculation leading to the prove of the Witten's conjecture in out particular
case.

\subsection{}
At the end, we have to acknowledge that we don't understand completely the original
formulation of the Witten's conjecture~\cite{w} and the more rigorous formulation
in~\cite{jkv} as well. But for our arguments in is enough to know the factorization
property of the Witten's top Chern class (see definitions below) and that everything is
already proved in genus zero.

This paper could be consider as a sequel to~\cite{s}. Though we have tried to make this
paper be self-closed, it could be useful to look through~\cite{s} to learn some
geometrical ideas hidden behind our calculations.

\section{The Witten's conjecture}\label{witcon}

\subsection{}
In this section we briefly formulate the Witten's conjecture. Note that
we consider here only the special case of the conjecture related to
the Boussinesq hierarchy. Note also that we doesn't follow the
initial Witten's definitions, but we use an axiomatic approach
from~\cite{jkv}.

\subsection{}
Consider the moduli space $\oM_{g,s}\ni(C,x_1,\dots,x_s)$.
Label each marked point $x_i$ by the number $m_i\in\{0,1\}$.

By $K$ denote the canonical line bundle of $C$. Consider the line bundle
$S=K\otimes(\otimes_{i=1}^s\O(x_i)^{-m_i})$ over $C$. If $2g-2-\sum_{i=1}^sm_i$
is divisible by $3$, then there are $3^{2g}$ isomorphism classes of line bundles
$\T$ such that $\T^{\otimes 3}\cong S$.

The choice of an isomorphism class of $\T$ determines a cover $\M'_{g,s}$ of
$\M_{g,s}$. To extend it to a covering of $\oM_{g,s}$ we have to discuss the
behavior of $\T$ near a double point.

\subsection{}
There are $3$ possible cases of behavior of $\T$ near the double point.
Two cases we call ''good'', and the third case is ''bad''.
For our purposes we need a description only of the good cases.

Let $C$ be a singular curve with one double point. By $\pi\colon C_0\to C$
denote its normalization. The preimage of the double point consists of two
points, say $x'$ and $x''$. Two good cases look like follows:
$\T\cong\pi_*\T'$, where $\T'$ is
a locally free sheaf on $C_0$ with a natural isomorphism either ${\T'}^{\otimes r}\cong
K\otimes(\otimes_{i=1}^s\O(x_i)^{-m_i})
\otimes\O(x')^{-1}$, or
${\T'}^{\otimes r}\cong
K\otimes(\otimes_{i=1}^s\O(x_i)^{-m_i})
\otimes\O(x'')^{-1}$.

\subsubsection{} Now we are ready to define the Witten's top Chern class
$c_D(\mathcal{V})\in H^{2D}(\oM'_{g,s})$,
which plays the main role in the Witten's cojecture.
The number $D$ is equal to $(g-1)/3+\sum_{i=1}^sm_i/3$.

If $g=0$, then $c_D(\mathcal{V})$ is the top Chern class of the bundle over $\oM'_{g,s}$,
which fiber at a moduli point $(C,x_1,\dots,x_s,\T)$ is equal to $H^1(C,\T)$.
In higher genera the Witten's top Chern class is just a cohomology class satisfying
the following properties.

First, consider a component of the boundary of $\oM'_{g,s}$ consisting of curves with
one double point, where $\T$ is defined by a ''good'' case.
If this component of the boundary consists of two-component curves,
then $c_D(\V)=c_{D_1}(\V_1)\cdot c_{D_2}(\V_2)$, where
$c_{D_i}(V_i)$ are the Witten's top Chern classes on the corresponding moduli spaces
of the components of normalizations.
If this component of the boundary consists of one-component self-intersecting curves, then
$c_D(\V)=c_D(\V')$, where $c_D(\V')$ is defined on the normalization of these
curves.

Second, consider a component of the boundary consisting of curves with one double point,
where we have the ''bad'' case of the behavior of $\T$. In this case we require that the
rescriction of $c_D(\V)$ to this component of the boundary vanishes.

\subsection{}
These properties uniquely determine the intersection numbers, which we consider
in this paper. As far as we understand, the existence of the Witten's top Chern class
is proved in~\cite{pv}.

\subsection{}
Let us label each marked point $x_i$ by an integer $n_i\geq 0$.
By $\l\tau_{n_1,m_1}\dots\tau_{n_s,m_s}\r_g$
denote the intersection number
\begin{equation}
\l\tau_{n_1,m_1}\dots\tau_{n_s,m_s}\r_g :=
\frac{1}{3^g}\int_{\oM'_{g,s}}\prod_{i=1}^s\psi_i^{n_i}\cdot
c_D(\mathcal{V}).
\end{equation}

Here $\psi_i$ is the first Chern class of the line bundle over $\oM'_{g,s}$,
whose fiber at a moduli point $(C,x_1,\dots,x_s,\T)$ is equal to $T^*_{x_i}C$.
Of course, the number $\l\tau_{n_1,m_1}\dots\tau_{n_s,m_s}\r_g$ is not zero only if
$3g-3+s=\sum_{i=1}^sn_i+D$.

\subsection{}
The special case of the conjecture from~\cite{w} is that the generating function
for the numbers $\l\tau_{n_1,m_1}\dots\tau_{n_s,m_s}\r_g$ is the string solution
to the Boussinesq hierarchy.

In particular, for the numbers $\l\tau_{n,m}\tau_{0,1}^k\tau_{0,0}^l\r_3$
the conjecture means that
\begin{equation}
\l\tau_{n,m}\tau_{0,1}^k\tau_{0,0}^l\r_3=
\frac{1}{3!12^3}\l\tau_{n-6,m}\tau_{0,1}^{k+3}\tau_{0,0}^l\r_0.
\end{equation}
(e.g., see relation for the Boussinesq hierarchy in~\cite{s})

This is just what we prove in this paper.

\section{Relations for intersection numbers}\label{relint}

\subsection{}
First of all, let us define some natural subvarieties of the moduli
spaces of curves.

Consider the moduli space $\oM_{g,s}'$. We assume that the covering
$\pi\colon\oM_{g,s}'\to\oM_{g,s}$ is defined by the collections of numbers
$m_1,\dots,m_s$. Of course, $m_i$ corresponds to the $i$-th marked
point, and each $m_i$ equals either to $0$ or to $1$.

Let $t<s$ be a positive integer. Consider a collection of positive integers
$a_1,\dots,a_t$. By $F(a_1,\dots,a_t)$ denote the subvariety of the open
moduli space $\M_{g,s}$ consisting of curves $(C,x_1,\dots,x_s)$ such that
$(-\sum_{i=1}^ta_i)x_1+a_1x_2+\dots+a_tx_{t+1}$ is a divisor of a meromorphic
function. By $F'(a_1,\dots,a_t)$ denote the closure of
$\pi^{-1}(F(a_1,\dots,a_t))$ in $\oM_{g,s}'$.

\subsection{}
Now let us introduce notations for the integrals over $F'(a_1,\dots,a_t)$.

Consider the numbers $m_{t+2},\dots,m_s$. Let $q$ of these numbers be equal
to $0$, and $p$ of these numbers be equal to $1$. Of course, $p+q=s-t-1$.
By $U(g,n,m_1,p,q|\prod_{i=1}^t\eta_{m_{i+1},a_i})$
denote the integral of $\psi_1^n\cdot c_D(\V)$ over the subspace $F'(a_1,\dots,a_t)$:

\begin{equation}
U(g,n,m_1,p,q|\prod_{i=1}^t\eta_{m_{i+1},a_i}) :=
\int_{F'(a_1,\dots,a_t)}\psi_1^n \cdot c_D(\V).
\end{equation}

\subsection{}
These numbers play the main role in our calculations below. In particular,
there is an expression of $\l\tau_{n,m}\tau_{0,1}^k\tau_{0,0}^l\r_3$ via
such numbers.

\begin{theorem}\label{theorem1}
\begin{multline}\label{formula1}
3!\l\tau_{n,m}\tau_{0,1}^k\tau_{0,0}^l\r_3=
U(3,n+1,m,k,l|\eta_{0,1}^4)
-3 U(3,n,m,k,l|\eta_{0,1}^3) \\
+3 U(3,n-1,m,k,l|\eta_{0,1}^2).
\end{multline}
\end{theorem}

\subsection{}
There is a recursion relation for the numbers
$U(g,n,m,p,q|\prod_{i=1}^t\eta_{m_i,a_i})$. Let us write it down
in the most general form.

\begin{theorem}\label{theorem2}
\begin{multline}\label{formula2}
(\sum_{i=1}^ta_i)(2g+t-1)
U(g,n,m,p,q|\prod_{i=1}^t\eta_{m_i,a_i}) = \\
\sum_{j:\, m_j=0}\sum_{b_1+b_2=a_j}
b_1 b_2 U(g-1,n-1,m,p,q|\eta_{1,b_1}\eta_{0,b_2}\prod_{i\not=j}\eta_{m_i,a_i}) + \\
\sum_{j:\, m_j=1}\sum_{b_1+b_2=a_j}
\frac{b_1 b_2}{2} U(g-1,n-1,m,p,q|\eta_{1,b_1}\eta_{1,b_2}\prod_{i\not=j}\eta_{m_i,a_i}) + \\
\sum_{j:\, m_j=1}\sum_{b_1+b_2=a_j}
\frac{p\, b_1 b_2}{6}
U(g-1,n-1,m,p-1,q|\eta_{0,b_1}\eta_{0,b_2}\prod_{i\not=j}\eta_{m_i,a_i}) + \\
\sum_{j:\, m_j=1}\sum_{b_1+b_2+b_3=a_j}
\frac{b_1 b_2 b_3}{9}
U(g-2,n-1,m,p,q|\eta_{0,b_1}\eta_{0,b_2}\eta_{0,b_3}\prod_{i\not=j}\eta_{m_i,a_i}) + \\
\sum_{j<k:\, m_j=m_k=0}
(a_j+a_k) U(g,n-1,m,p,q|\eta_{0,a_j+a_k}\prod_{i\not=j,k}\eta_{m_i,a_i}) + \\
\sum_{j:\, m_j=0}\sum_{k:\, m_k=1}
(a_j+a_k) U(g,n-1,m,p,q|\eta_{1,a_j+a_k}\prod_{i\not=j,k}\eta_{m_i,a_i}) + \\
\sum_{j<k:\, m_j=m_k=1}
\frac{p (a_j+a_k)}{3} U(g,n-1,m,p-1,q|\eta_{0,a_j+a_k}\prod_{i\not=j,k}\eta_{m_i,a_i}) + \\
\sum_{j<k:\, m_j=m_k=1}\sum_{b_1+b_2=a_j+a_k}
\frac{b_1 b_2}{3} U(g-1,n-1,m,p,q|\eta_{0,b_1}\eta_{0,b_2}\prod_{i\not=j,k}\eta_{m_i,a_i}) + \\
\sum_{j<k<l:\, m_j=m_k=m_l=1}
\frac{2(a_j+a_k+a_l)}{3} U(g,n-1,m,p,q|\eta_{0,a_j+a_k+a_l}\prod_{i\not=j,k,l}\eta_{m_i,a_i}).
\end{multline}
\end{theorem}

\subsection{}
In fact, both these theorems are very simple corollaries of the technique
developed in~\cite{i}. This is discussed in details in~\cite{s}. Moreover,
if we put $p=q=0$, then both these theorems appear to be special cases of
the similiar theorems in~\cite{s}.

Since both these theorems could be proved just by the same argument as given
in~\cite{s}, we skip the proofs.

\section{Calculations}\label{calcul}

\subsection{}
In this section we prove the following theorem:

\begin{theorem}\label{theorem3}
\begin{equation}
\l\tau_{n,m}\tau_{0,1}^k\tau_{0,0}^l\r_3=
\frac{1}{3!12^3}\l\tau_{n-6,m}\tau_{0,1}^{k+3}\tau_{0,0}^l\r_0.
\end{equation}
\end{theorem}

\subsection{}
To prove this theorem we need two lemmas.

\begin{lemma}
\begin{equation}
\l\tau_{n,m}\tau_{0,1}^{k}\tau_{0,0}^l\r_0=
\l\tau_{n+1,m}\tau_{0,1}^k\tau_{0,0}^{l+1}\r_0
\end{equation}
\end{lemma}

This is just the special case of the string equation (see~\cite{w}).

\begin{lemma}
\begin{equation}
U(0,n,m,p,q|\prod_{i=1}^t\eta_{m_{i},a_i})=
\l\tau_{n,m}\tau_{0,1}^p\tau_{0,0}^q\prod_{i=1}^t\tau_{0,m_{i}}\r_0.
\end{equation}
\end{lemma}

This lemma immediately follows from our definitions.

\subsection{}
Applying many times theorem~\ref{theorem2} and these two lemmas we obtain

\begin{multline}
U(3,n+1,m,k,l|\eta_{0,1}^4)=
\frac{317}{1548288}\l\tau_{n-6,m}\tau_{0,1}^{k+3}\tau_{0,0}^l\r_0+\\
\left(\frac{11791}{23224320}+\frac{769k}{1548288}\right)
\l\tau_{n-7,m}\tau_{0,1}^{k}\tau_{0,0}^{l+1}\r_0;
\end{multline}

\begin{multline}
U(3,n,m,k,l|\eta_{0,1}^3)=
\frac{1}{30618}\l\tau_{n-6,m}\tau_{0,1}^{k+3}\tau_{0,0}^l\r_0+\\
\left(\frac{919}{7348320}+\frac{223k}{1837080}\right)
\l\tau_{n-7,m}\tau_{0,1}^{k}\tau_{0,0}^{l+1}\r_0;
\end{multline}

\begin{equation}
U(3,n-1,m,k,l|\eta_{0,1}^2)=
\frac{k+1}{120960}\l\tau_{n-7,m}\tau_{0,1}^{k}\tau_{0,0}^{l+1}\r_0.
\end{equation}

\subsection{}
It follows from these relations and theorem~\ref{theorem1} that
\begin{multline}\label{final}
3!\l\tau_{n,m}\tau_{0,1}^k\tau_{0,0}^l\r_3=
\frac{446}{41803776}\l\tau_{n-6,m}\tau_{0,1}^{k+3}\tau_{0,0}^l\r_0+\\
\frac{19729(k+1)}{125411328}
\l\tau_{n-7,m}\tau_{0,1}^{k}\tau_{0,0}^{l+1}\r_0;
\end{multline}

\subsection{}
The standard calculations in genus zero (e.g., see~\cite{w}) imply the following lemma:

\begin{lemma}\label{lemma3}
\begin{equation}
\l\tau_{n,m}\tau_{0,1}^{k+3}\tau_{0,0}^l\r_0=
\frac{k+1}{3}\l\tau_{n-1,m}\tau_{0,1}^k\tau_{0,0}^{l+1}\r_0
\end{equation}
\end{lemma}

Combining equation~\ref{final} and lemma~\ref{lemma3}, we get
\begin{equation}
3!\l\tau_{n,m}\tau_{0,1}^k\tau_{0,0}^l\r_3=
\frac{1}{1728}\l\tau_{n-6,m}\tau_{0,1}^{k+3}\tau_{0,0}^l\r_0.
\end{equation}

This concludes the prove.

\bigskip

e-mail: shadrin@mccme.ru

\end{document}